\documentclass[preprint,12pt, a4paper]{elsarticle}

\usepackage{amssymb}
\usepackage{hyperref}
\usepackage{amsmath}
\usepackage{listings}
\usepackage{xcolor}

\definecolor{codegreen}{rgb}{0,0.6,0}
\definecolor{codegray}{rgb}{0.5,0.5,0.5}
\definecolor{codepurple}{rgb}{0.58,0,0.82}
\definecolor{backcolour}{rgb}{0.95,0.95,0.92}

\lstdefinestyle{mystyle}{
    backgroundcolor=\color{backcolour},   
    commentstyle=\color{codegreen},
    keywordstyle=\color{magenta},
    numberstyle=\tiny\color{codegray},
    stringstyle=\color{codepurple},
    basicstyle=\ttfamily\footnotesize,
    breakatwhitespace=false,         
    breaklines=true,                 
    captionpos=b,                    
    keepspaces=true,                 
    numbers=left,                    
    numbersep=5pt,                  
    showspaces=false,                
    showstringspaces=false,
    showtabs=false,                  
    tabsize=2
}

\lstdefinelanguage{julia}
{
  keywordsprefix=\@,
  morekeywords={
    exit,whos,edit,load,is,isa,isequal,typeof,tuple,ntuple,uid,hash,finalizer,convert,promote,
    subtype,typemin,typemax,realmin,realmax,sizeof,eps,promote_type,method_exists,applicable,
    invoke,dlopen,dlsym,system,error,throw,assert,new,Inf,Nan,pi,im,begin,while,for,in,return,
    break,continue,macro,quote,let,if,elseif,else,try,catch,end,bitstype,ccall,do,using,module,
    import,export,importall,baremodule,immutable,local,global,const,Bool,Int,Int8,Int16,Int32,
    Int64,Uint,Uint8,Uint16,Uint32,Uint64,Float32,Float64,Complex64,Complex128,Any,Nothing,None,
    function,type,typealias,abstract, DataFrame, get,round,sum,find_solution,plot,plot!,
    find_search_interval
  },
  sensitive=true,
  morecomment=[l]{\#},
  morestring=[b]',
  morestring=[b]" 
}

\journal{SoftwareX}

\begin{document}
\renewcommand{\labelenumii}{\arabic{enumi}.\arabic{enumii}}

\begin{frontmatter}
 


\title{TulipaProfileFitting.jl: A Julia package for fitting renewable energy time series profiles}


\author[label1]{Diego A. Tejada-Arango}
\author[label2]{Germán Morales-España}
\author[label3]{Abel S. Siqueira}
\author[label4]{Özge Özdemir}
\address[label1]{TNO, Motion Building, Radarweg 60, 1043 NT Amsterdam, diego.tejadaarango@tno.nl}
\address[label2]{TNO, Motion Building, Radarweg 60, 1043 NT Amsterdam, german.morales@tno.nl}
\address[label3]{The Netherlands eScience Center, Science Park 402, 1098 XH Amsterdam, abel.siqueira@esciencecenter.nl}
\address[label4]{PBL, Bezuidenhoutseweg 30, 2594 AV Den Haag, ozge.ozdemir@pbl.nl}

\begin{abstract}
This paper introduces the TulipaProfileFitting.jl package, a tool developed in Julia for generating renewable energy profiles that fit a specified capacity factor. It addresses the limitations of naive methods in adjusting existing profiles to match improved technology efficiency, particularly in scenarios lacking detailed weather data or technology specifications. By formulating the problem mathematically, the package provides a computationally efficient solution for creating realistic renewable energy profiles based on existing data. It ensures that the adjusted profiles realistically reflect the improvements in technology efficiency, making it an essential tool for energy modellers in analyzing $CO_2$-neutral energy systems.
\end{abstract}


\end{frontmatter}


\section*{Metadata}
\label{}

\begin{table}[!h]
\begin{tabular}{|l|p{6.5cm}|p{6.5cm}|}
\hline
\textbf{Nr.} & \textbf{Code metadata description} & \textbf{Please fill in this column} \\
\hline
C1 & Current code version & v0.3.0 \\
\hline
C2 & Permanent link to code/repository used for this code version & \url{https://github.com/TulipaEnergy/TulipaProfileFitting.jl} \\
\hline
C3  & Permanent link to Reproducible Capsule & \\
\hline
C4 & Legal Code License   & Apache License 2.0 \\
\hline
C5 & Code versioning system used & git \\
\hline
C6 & Software code languages, tools, and services used & Julia \\
\hline
C7 & Compilation requirements, operating environments \& dependencies & Julia\\
\hline
C8 & If available Link to developer documentation/manual & \url{https://tulipaenergy.github.io/TulipaProfileFitting.jl/dev/} \\
\hline
C9 & Support email for questions & diego.tejadaarango@tno.nl\\
\hline
\end{tabular}
\caption{Code metadata}
\label{codeMetadata} 
\end{table}

\section{Motivation and significance}
Energy system models, like SpineOpt \cite{Ihlemann2022} and PyPSA \cite{Horsch2018}, are widely used to analyze different scenarios and pathways for transitioning towards $CO_2$-neutral energy systems. A significant part of this transition depends on weather-dependent renewable energy sources, such as wind and solar. Therefore, renewable energy profiles, representing the available production of these resources, are crucial inputs for these models. Practitioners use raw weather data, such as wind speed or solar radiation, to determine these profiles \cite{Koivisto2020}. Then, they determine the availability profile according to the technology characteristics, such as turbine type and height for wind, and tilt and tracking system for solar \cite{Staffell2016}. Capacity factors, the mean values of these profiles in a year, indicate the efficiency of the technology in producing energy. Renewable energy resources with higher capacity factors are more efficient to take advantage of the weather resources. With improvements in renewable technology, capacity factors are expected to increase by nearly a third to 31.3\% \cite{Staffell2016}.

As an energy modeller, if you have the raw weather data and know the technical details of the new technology for creating a new profile, you can use the methods developed in \cite{Koivisto2020} and \cite{Staffell2016} and then use the results in your energy model for analysis. However, what if you do not have the weather data, are unaware of the specific characteristics of the new technology, or need a fast but good approximation of the improvement in the capacity factor? A naive approach is to multiply an existing renewable profile by a factor that scales up the data until the target capacity factor is achieved. However, this approach is flawed since it could result in unrealistic values. For example, if your profile has values equal to 1 p.u. and you consider an improvement of 5\%, then you will end up with a profile value of 1.05 p.u. Even if you end up with values below 1 p.u. using the naive approach, the higher values still receive a greater increase than the lower ones, which is not realistic since technology efficiency improvements are typically allocated in the middle values of the profiles.

To overcome this problem, the authors in \cite{Ozdemir2020} used an optimization approach to obtain new profiles without the drawbacks of the naive approach. However, our findings suggest that it is possible to achieve similar results without formulating an optimization problem. Therefore, we introduce the \textit{TulipaProfileFitting.jl} package, developed in Julia \cite{Julia-2017}, which is designed to create new renewable energy profiles from existing ones fitting them to a target capacity factor. We define the problem as a mathematical equation and reformulate it to find a coefficient fitting the existing renewable profile with a target capacity factor. This approach simplifies the problem, making it more computationally efficient and relevant for real-world applications. The package includes a robust algorithm that identifies the most appropriate values for a new profile, ensuring that the average energy output aligns with the target capacity factor chosen by the user.

\section{Software description}
\textit{TulipaProfileFitting.jl} is a package based on the mathematical equation (\ref{fig:eq-fx}). The equation considers a list of $m$ numbers representing potential renewable energy production values ranging from $0$ to $1$. The primary challenge is finding a power $x$ that satisfies the condition where the mean of these values, each raised to $x$, equals a target mean $\mu \in (0, 1)$.

\begin{equation}\label{fig:eq-fx}
  \frac{1}{m}\sum_{i=1}^m p_i^x = \mu  
\end{equation}

Using equation (\ref{fig:eq-fx}), we can analyze the following questions to solve the problem:

\begin{enumerate}
    \item[Q1.] Is it possible to find such a value of $x$?
    \item[Q2.] If such $x$ exists, how do we find it?
    \item[Q3.] What should be done if it is not possible to find such a value of $x$?
\end{enumerate}

Both $p_i$ and $x$ can be zero, resulting in an undefined value of $p_i^x$. Some programming languages, including Julia, define $0^0$ as $1$. However, this definition does not make sense for the fitting profile application, as the $p_i$ values are fixed while $x$ is variable. As a result, we define the function $\sigma$ as $\sigma:[0,1]\times[0,\infty)$ to avoid this issue.

\begin{equation}\label{fig:eq-sigma}
    \sigma(p, x) = \left\{\begin{array}{ll}p^x, \ & \text{if } p > 0, \\ 0, \ & \text{otherwise}.\end{array}\right.
\end{equation}

Equation (\ref{fig:eq-sigma}) This allows $x \mapsto \sigma(p, x)$ to be continuous for any value $p \in [0, 1]$.

Therefore, the problem can be redefined as finding $x$ such that:

$$\frac{1}{m}\sum_{i=1}^m \sigma(p_i, x) = \mu.$$

However, to avoid this notation, we can assume, without loss of generality, that there are $r$ non-zero $p_i$, i.e., $p_1 \geq \cdots \geq p_r > 0$, and $p_{r+1} = \cdots = p_m = 0$.

This simplifies the problem of finding $x$ such that:

$$\frac{1}{m} \sum_{i=1}^r p_i^x = \mu.$$

Let's define $S(x) = \frac{1}{m}\sum_{i=1}^r p_i^x$ to help us with the notation. To visualize the problem, we are going to use four sets of possible profile values:

\begin{itemize}
    \item A set with a 0 and a 1
    \item A set with a lot of zeros and ones (25\% each)
    \item A set that only contains numbers from 0 to $\alpha < 0.5$.
    \item A set that only contains numbers from $(1 - \alpha)$ to 1.
\end{itemize}

\begin{figure}
    \centering
    \includegraphics[width=1.0\linewidth]{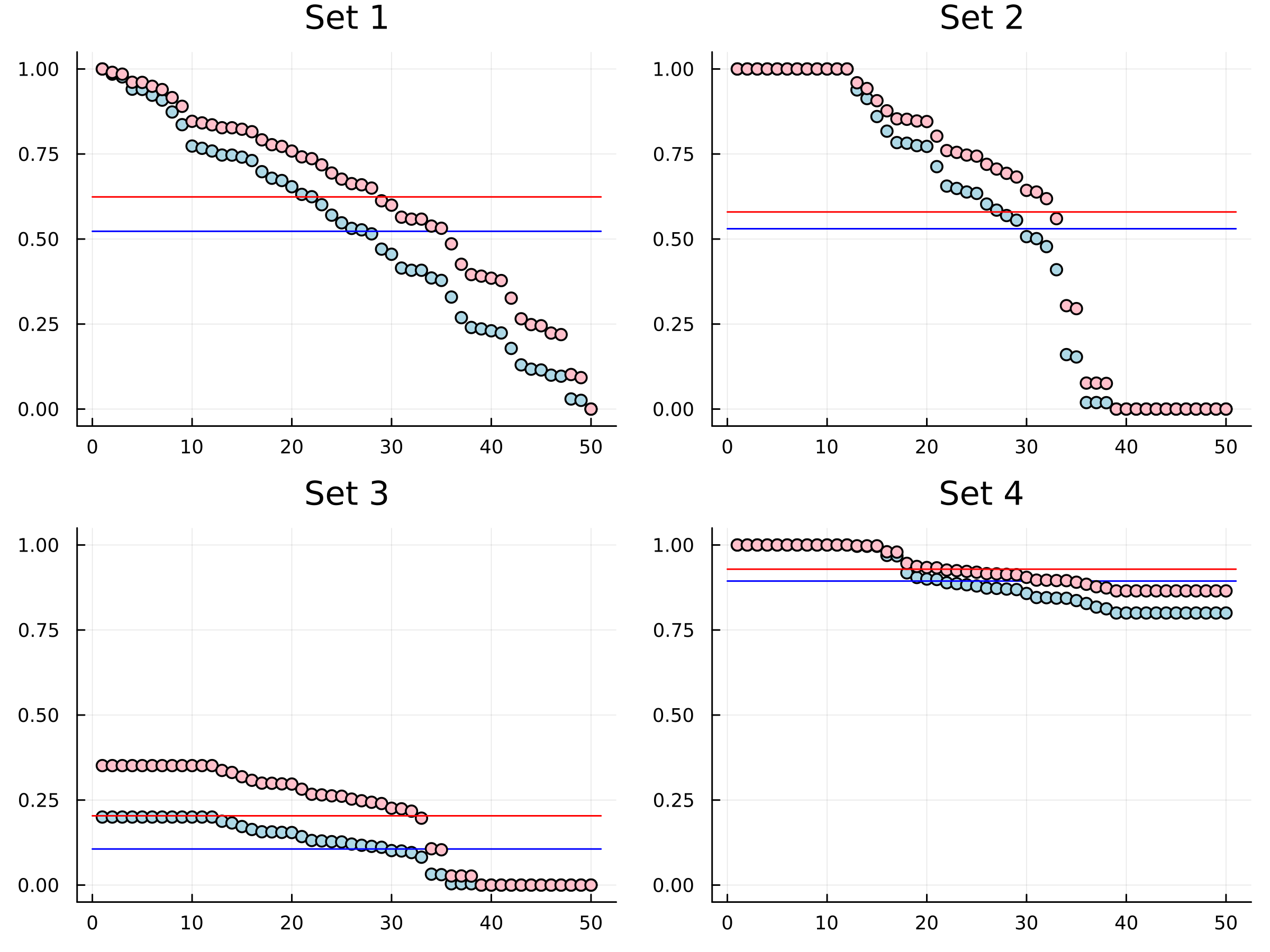}
    \caption{Example of sets of values}
    \label{fig:ex-sets}
\end{figure}

First, let's visualize the effect of $p_i^x$ on these sets. Figure (\ref{fig:ex-sets}) displays the original values in \textit{lightblue}, sorted from the maximum to the minimum. The mean value of the original data is indicated in \textit{blue}. In addition, the value for $x=0.5$ is indicated in \textit{pink}, and the new mean is depicted in \textit{red}.

\begin{figure}
    \centering
    \includegraphics[width=0.75\linewidth]{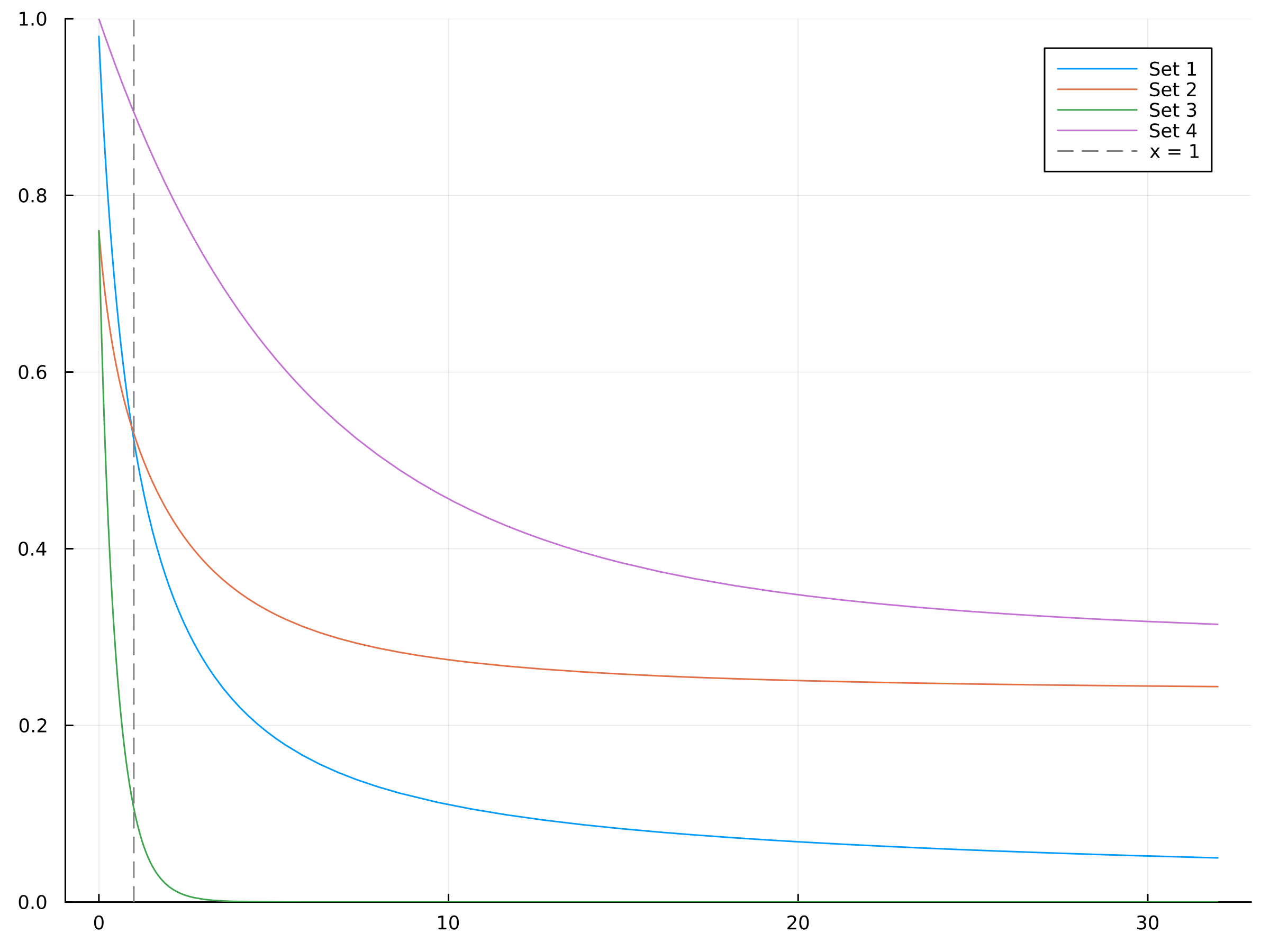}
    \caption{Function $S$ for each set of values}
    \label{fig:ex-sets-functions}
\end{figure}

Furthermore, Figure (\ref{fig:ex-sets-functions}) shows the function $S$ for each set.

\textit{Answer to Q1}: Since $p^0 = 1$ for positive $p$, then $S(0) = \frac{r}{m}$. $S$ is non-increasing since

$$S'(x) = \frac{1}{m} \sum_{i=1}^r p_i^x \ln p_i \leq 0.$$

And assuming that there are $n$ values such that $p_i = 1$, then:

$$\lim_{x \to \infty} S(x) = \frac{n}{m}.$$

This means that there is a solution to the problem if $\dfrac{n}{m} < \mu \leq \dfrac{r}{m}$.

\textit{Answer to Q2}: Assuming $S$ decreasing and $\mu$ in range, we can solve the problem by looking for an interval $[a, b]$ such that $S(x) - \mu$ changes sign.
This can be done by creating an increasing sequence of points $v_1, v_2, \dots$, such that $v_1 = 0$, and $v_{i+1} > v_i$ with $\lim_{i \to \infty} v_i = \infty$.
For instance, the sequence $0, 1, 2, 4, 8, \dots$. Since $S$ is decreasing, and $\mu$ is in range, then either $S(v_i) = \mu$ for some $v_i$, or $S(x) - \mu$ will change sign in some interval $[v_i, v_{i+1}]$. Given the interval, we perform a bisection search or similar.

\textit{Answer to Q3}: If $\mu > \dfrac{r}{m}$, then $x = 0$ will yield $S(0) = \dfrac{r}{m}$, which is the closest to $\mu$. Alternatively, if $\mu \leq \dfrac{n}{m}$, then $S(x) \to \dfrac{n}{m}$ as $x \to \infty$. In this case, we select a reasonably large value for $x$.

\begin{figure}
    \centering
    \includegraphics[width=1.0\linewidth]{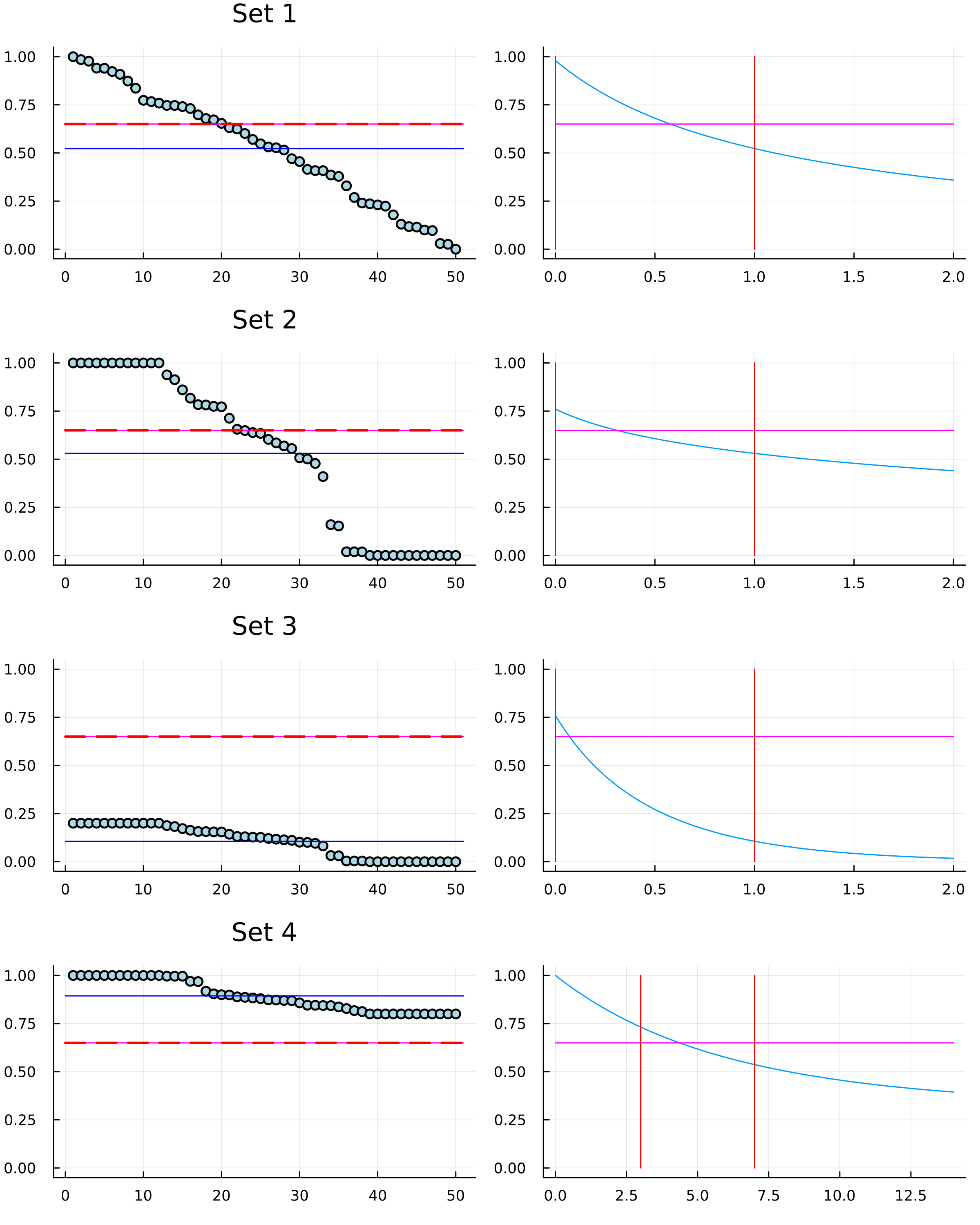}
    \caption{Results for four sets of values}
    \label{fig:ex-sets-results}
\end{figure}

Figure (\ref{fig:ex-sets-results}) shows the results of each set with a target mean of $\mu = 0.65$. The plots on the left-hand side highlight the target value in \textit{red}. Meanwhile, the plots on the right-hand side depict the interval $\dfrac{n}{m} < \mu \leq \dfrac{r}{m}$. We can observe that the function $S$ of each set of values intersects with the target value.

\subsection{Software architecture}
The package was developed in Julia using the Root.jl \cite{Roots.jl} and Statistics.jl. It is registered in the Julia and can be installed using the following commands in the Julia REPL.  

\begin{lstlisting}[language=julia]
using Pkg
Pkg.add("TulipaProfileFitting")
\end{lstlisting}

\subsection{Software functionalities}
The package has the following functions or methods:

\begin{itemize}
    \item $TulipaProfileFitting.find\_search\_interval$: Returns an interval such that $f(a) \cdot f(b) > 0$. It could be 0 for either endpoint, but it is not positive, ensuring that there is a root in $[a, b]$.
    \item $TulipaProfileFitting.find\_solution$: Finds a points such that $S(x) = \mu$, if possible, where $S(x) = \frac{1}{|P|} \sum_{p \in P: p > 0} p^x$. If not possible, return either 0 or 1000, depending on what is most appropriate.
    \item $TulipaProfileFitting.validate_profile$: Validates if the values of profile P are within 0 and 1. If not, it throws an error message and stops the calculation.
\end{itemize}

\section{Illustrative examples}
\textit{TulipaProfileFitting.jl} is a package designed for fitting power production curves of renewable sources, such as wind and solar. In this particular example, we have generated a time series for a wind power plant's power production, using a method that was developed in reference \cite{Staffell2016}. The code below illustrates how to use the package and plot its results.

\begin{lstlisting}[language=julia]
# Load packages
using TulipaProfileFitting
using CSV
using Plots
using DataFrames
using HTTP

plotly()

# File location
file_url = "https://raw.githubusercontent.com/TulipaEnergy/TulipaProfileFitting.jl/main/docs/src/files/wind_power_profile.csv"

# Read file
df = DataFrame(CSV.File(HTTP.get(file_url).body, header=4))

# Get the profile from the dataframe
profile_values = df.electricity

# Current capacity factor (e.g., mean value)
current_cp = round(sum(profile_values)/8760;digits=2)

# New capacity factor as definition
target_cp = 0.6

# Obtain the coefficient that fits the values to the target.
coefficient = find_solution(profile_values, target_cp)

# Determine the fitted profile
fitted_profile = profile_values.^coefficient

# Plot chronological values
plot(profile_values, label="profile")
plot!(fitted_profile, label="fitted")

# Plot sorted values
plot(sort(profile_values,rev=true), label="profile")
plot!(sort(fitted_profile,rev=true), label="fitted")
\end{lstlisting}

\begin{figure}
    \centering
    \includegraphics[width=1.0\linewidth]{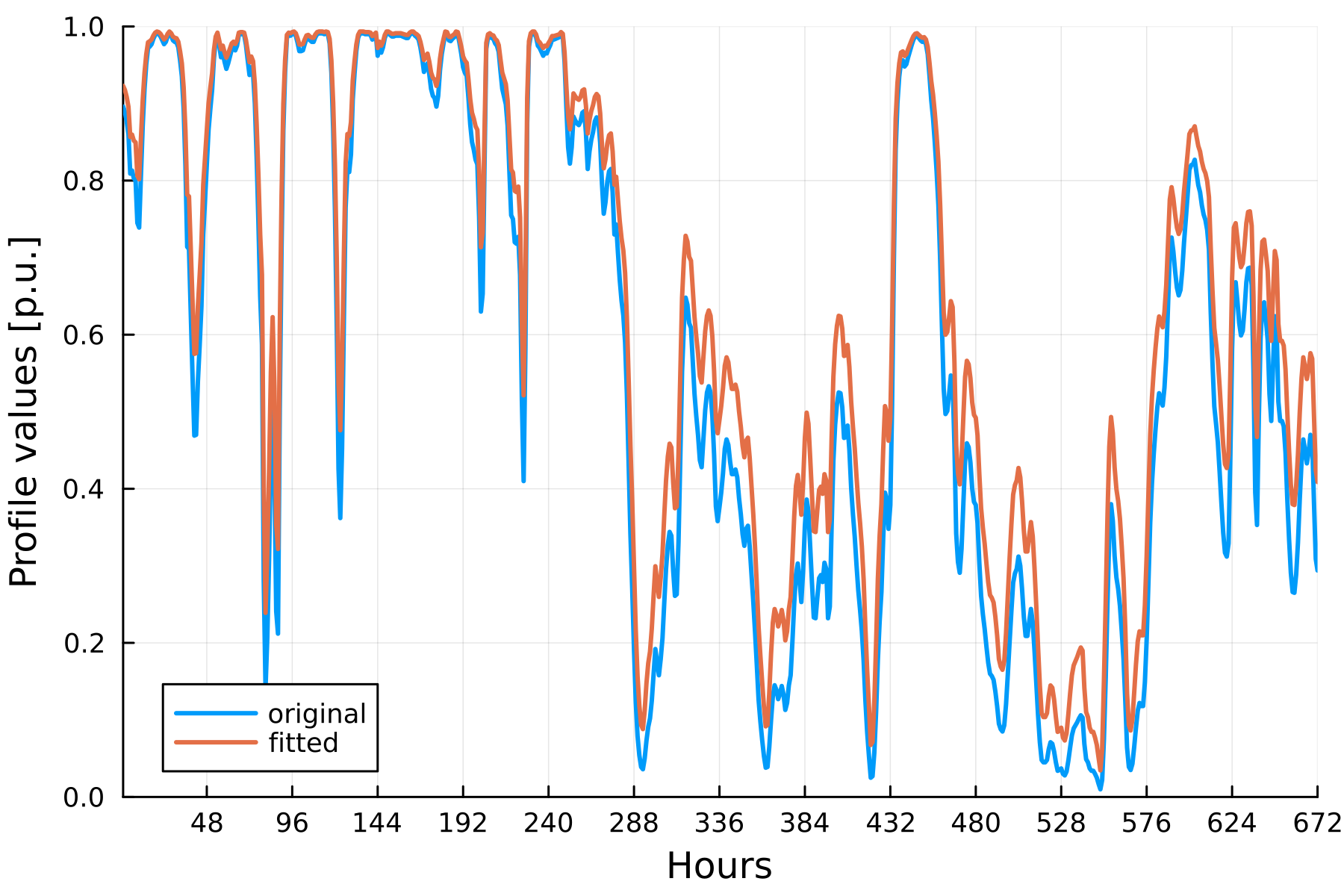}
    \caption{Example of the hourly profile for one week}
    \label{fig:ex-hourly}
\end{figure}

\begin{figure}
    \centering
    \includegraphics[width=1.0\linewidth]{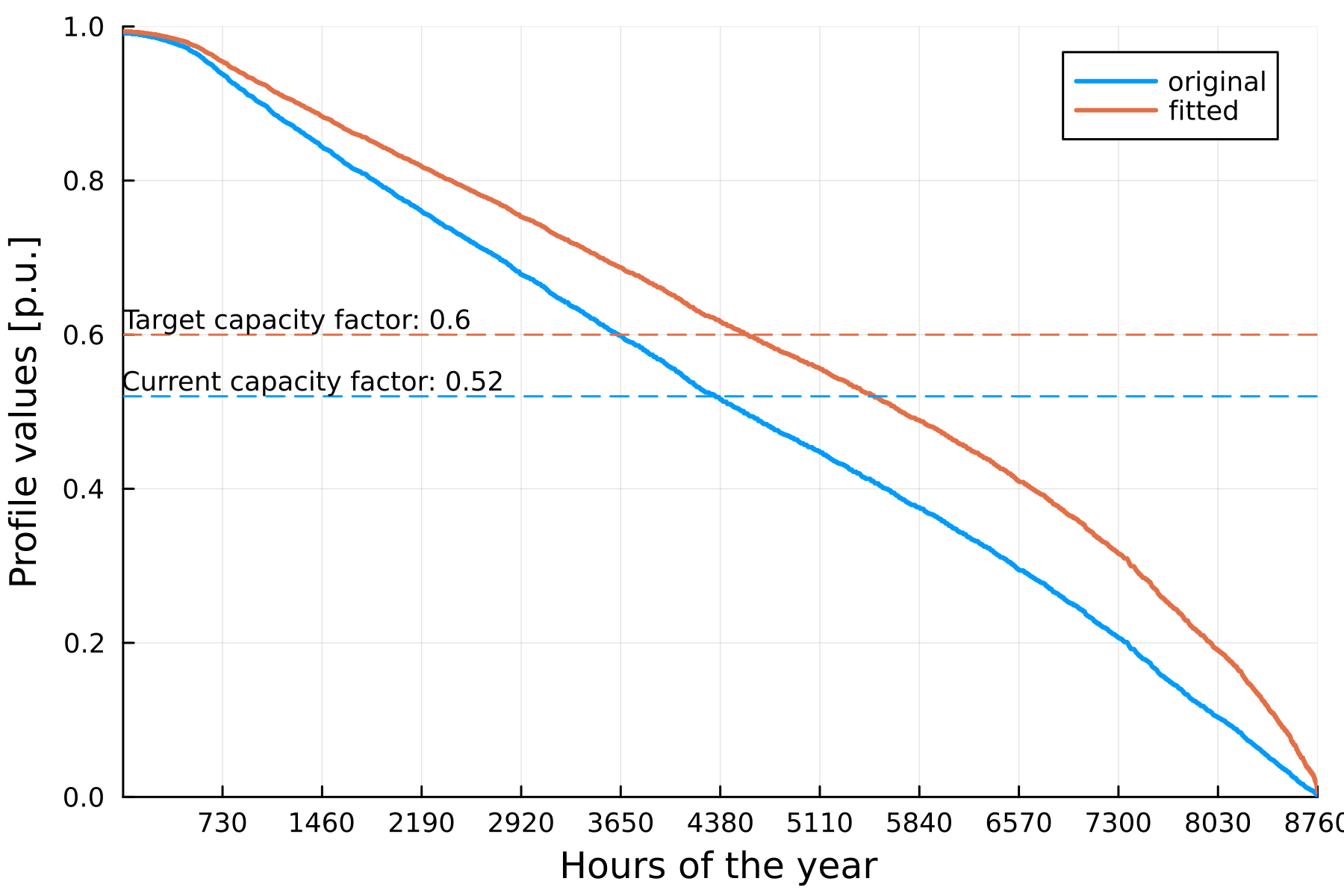}
    \caption{Sorted values for the profile in the year}
    \label{fig:ex-sorted}
\end{figure}

The example results are shown in Figures \ref{fig:ex-hourly} and \ref{fig:ex-sorted}. It can be noticed that the fitted curve mainly affects the middle values within the entire range, while the values closer to zero or one remain relatively unchanged. This outcome is expected since the package's main objective is to adjust the middle values, increasing or decreasing the capacity factor. Finally, Since the package uses well-known numerical methods, solution times are less than a second and can be parallelised for multiple profiles.

\section{Impact}
The package is designed to provide a quick and reliable method for generating new renewable energy profiles based on existing ones. It is particularly useful for complex scenarios that involve a combination of high and low-output renewable sources. The package serves as a practical tool for modellers who work with energy system models and wish to analyze multiple scenarios with improved capacity factors for renewable energy sources.

\section{Conclusions}
This package provides a mathematical representation and an intuitive solution for generating new renewable resource profiles that account for their improved capacity factor in a fast and reliable manner. Its rigorous mathematical foundation, combined with its practical application, positions it as a valuable tool for both researchers and practitioners in the energy modelling field. Using Julia as a programming language provides the package with the advantage of being in an open-source environment, allowing for further development and future research.

\section*{Acknowledgements}
This research received funding from the Dutch Research Council (NWO) under grant number ESI.2019.008. It reflects only the author’s views, and NWO is not liable for any use that may be made of the information contained therein.

The authors would like to express their gratitude to Professor Benjamin F. Hobbs from the Department of Environmental Health and Engineering at Johns Hopkins University for his insightful suggestion to use optimization techniques to fit renewable energy profiles. This innovative approach has formed the foundation of our research, and we are grateful for his valuable contribution.



\bibliographystyle{elsarticle-num} 
\bibliography{references}






\end{document}